\documentclass[12pt]{amsart}
\title{Schr\"odinger equation with linear potential and hitting times}
\author{Gerardo Hern\'andez-del-Valle}
\begin{document}
\maketitle
\begin{abstract}
In Hern\'andez-del-Valle (2010) the author studies the connection between Schr\"odinger's equation and first hitting densities of Brownian motion. Although the author is able to find solutions of a Schr\"odinger type pde he fails---except in some special cases---to construct a solution which satisfies the boundary on the space variable at $x=0$. In this paper we use an approach used in Bluman and Shtelen (1996) to find solutions which satisfy the pde and boundary condition when $t=0$.
\end{abstract}
Consider the Schr\"odinger equation
\begin{eqnarray}\label{schrodinger}
\frac{\partial u_1}{\partial t}(t,x)+\frac{\partial^2u_1}{\partial x^2}(t,x)-V_1(t,x)u_1(t,x)=0.
\end{eqnarray}

Given a linear operator $\mathcal{L}$, its adjoint $\mathcal{L}^*$ is defined by
$$\Phi\mathcal{L}u-u\mathcal{L}^*\Phi=\sum\limits_{i=1}^n D_if^i$$
where $x=(x_1,x_2,\dots,x_n)$, the total derivative operators $D_i=\partial/\partial x_i$, $i=1,2,\dots,n$, and
$\{f^i\}$ are bilinear expression in $u,\Phi$ and their derivatives. Consequently, if
\begin{equation}\label{adjoint}
\mathcal{L}^*\Phi=0
\end{equation}
then $\mathcal{L}u=0$ if and only if $\sum_{i=1}^nD_if_i=0$, i.e. a given linear partial differential equation
\begin{equation}\label{operator}
\mathcal{L}u=0
\end{equation}
is equivalent to the conservation law
\begin{equation}\label{conservation}
\sum\limits_{i=1}^nD_if^i=0
\end{equation}
for any $\Phi$ satisfying its adjoint equation (\ref{adjoint}).

We now specialize to the case when (\ref{operator}) is the Schr\"odinger equation (\ref{schrodinger}). Here the linear operator is
\begin{eqnarray*}
\mathcal{L}=\frac{\partial}{\partial t}+\frac{\partial^2}{\partial x^2}-V_1(t,x)
\end{eqnarray*}
its adjoint is given by
\begin{eqnarray*}
\mathcal{L}^*=-\frac{\partial}{\partial t}+\frac{\partial^2}{\partial x^2}-V_1(t,x)
\end{eqnarray*}
and (\ref{conservation}) becomes
\begin{eqnarray}\label{conservation1}
\frac{\partial}{\partial t}\left(\Phi u_1\right)+\frac{\partial}{\partial x}\left(\Phi\frac{\partial u_1}{\partial x}-\frac{\partial \Phi}{\partial x}u_1\right)=0.
\end{eqnarray}
The potential system corresponding to (\ref{conservation1}) is given by
\begin{eqnarray}
\label{v1}\frac{\partial v_1}{\partial x}&=& \Phi u_1\\
\label{v11}\frac{\partial v_1}{\partial t}&=&\frac{\partial \Phi}{\partial x}u_1-\Phi\frac{\partial u_1}{\partial x}
\end{eqnarray}
where $\Phi(t,x)$ is a solution of
\begin{eqnarray}\label{forward}
\mathcal{L}^*\Phi=-\frac{\partial\Phi}{\partial t}+\frac{\partial^2\Phi}{\partial x^2}-V_1(t,x)\Phi=0.
\end{eqnarray}

Note that if ($u_1(t,x)$, $v_1(t,x)$, $\Phi(t,x)$) solves (\ref{v1}), (\ref{v11}) and (\ref{forward}) then $u_1(t,x)$ solves Schr\"odinger equation (\ref{schrodinger}) and $v_1(t,x)$ solves
\begin{eqnarray}\label{v2}
\frac{\partial v_1}{\partial t}+\frac{\partial^2v_1}{\partial x^2}-\frac{2}{\Phi}\frac{\partial\Phi}{\partial x}\frac{\partial v_1}{\partial x}=0.
\end{eqnarray}
If $u_1(t,x)$ solves (\ref{schrodinger}) and $\Phi(t,x)$ solves (\ref{forward}), then one can $v_1(t,x)$ solving (\ref{v1}) and (\ref{v11}) i.e. for any $\Phi(t,x)$ satisfying (\ref{forward}) the point transformation
\begin{eqnarray}\label{w}
w=\frac{v_1}{\Phi}
\end{eqnarray}
maps (\ref{v2}) to
\begin{eqnarray}\label{pdew}
\frac{\partial w}{\partial t}+\frac{\partial^2w}{\partial x^2}-V_2(t,x) w=0
\end{eqnarray}
where the new potential $V_2(t,x)$ is given by
\begin{eqnarray}\label{V2}
V_2(t,x)=V_1(t,x)-2\frac{\partial^2}{\partial x^2}\log\Phi.
\end{eqnarray}

From equations (\ref{v1}), (\ref{v11}) and (\ref{w}) it follows that if $u_1(t,x)$ solves (\ref{schrodinger}) then
\begin{eqnarray}\label{solution}
w(t,x)=\frac{1}{\Phi(t,x)}\left[\int_k^xu_1(t,\xi)\Phi(t,\xi)d\xi+B_2(t)\right]
\end{eqnarray}
with $B_2(t)$ satisfying the condition
\begin{eqnarray}\label{B2}
\frac{dB_2}{dt}=\frac{\partial \Phi}{\partial x}(k,t)u_1(t,k)-\Phi(t,k)\frac{\partial u_1}{\partial x}(t,k)
\end{eqnarray}
for any constant $k$, solves the Schr\"odinger equation (\ref{w}).

In particular, let $V_1(t,x)=xf''(t)$ then (\ref{forward}) is
\begin{eqnarray*}
\frac{\partial \Phi}{\partial t}(t,x)+xf''(t)\Phi(t,x)=\frac{\partial^2\Phi}{\partial x^2}(t,x).
\end{eqnarray*}
which in turn admits solutions of the following form:
\begin{eqnarray*}
\Phi(t,x)=\exp\left\{\frac{1}{2}\int_0^t(f'(u))^2du-xf'(t)\right\}\omega\left(t,x-\int_0^tf'(s)ds\right)
\end{eqnarray*}
where $\omega$ in turn is a solution of 
\begin{eqnarray}\label{heat}
\omega_t=\omega_{xx}.
\end{eqnarray}
a particular solution to (\ref{heat}) is given by
\begin{eqnarray*}
\omega(t,x)=\exp\left\{-\frac{1}{2}\lambda^2t\pm \lambda x\right\}
\end{eqnarray*}
for some scalar $\lambda$, which alternatively implies that a solution to (\ref{w}) is given by:
\begin{eqnarray}\label{phi}
\Phi(t,x)=e^{\frac{1}{2}\int_0^t(f'(u))^2du-x[f'(t)\pm \lambda]-\frac{1}{2}\lambda^2t\pm \int_0^tf'(u)du}.
\end{eqnarray}
Next, note that the potential $V_2(t,x)=V_1(t,x)$ since 
$$\frac{\partial ^2}{\partial x^2}\log\Phi=0.$$
Hence (\ref{pdew}) becomes:
\begin{eqnarray}\label{pdew1}
-\frac{\partial w}{\partial t}(t,x)+xf''(t)w(t,x)=\frac{1}{2}\frac{\partial^2w}{\partial x^2}(t,x)
\end{eqnarray}
or from (\ref{solution}) and setting $k=0$
\begin{eqnarray*}
w(t,x)=\frac{1}{\Phi(t,x)}\left[\int_0^xu_1(t,\xi)\Phi(t,\xi)d\xi+B_2(t)\right],
\end{eqnarray*}
where $u_1$ solves Schr\"odinger's equation (\ref{schrodinger}), $\Phi$ is as in (\ref{phi}) and $B_2$ is defined in (\ref{B2}).

Solutions of $u_1$ are:
\begin{eqnarray*}
u_1(t,x)=\exp\left\{\frac{1}{2}\int_t^s(f'(u))^2du+xf'(t)\right\}\omega\left(s-t,x+\int_t^sf'(u)du\right).
\end{eqnarray*}
Then $w(0,x)$ becomes:
\begin{eqnarray*}
w(0,x)=\frac{1}{\Phi(0,x)}\int_0^xu_1(0,\xi)\Phi(0,\xi)d\xi
\end{eqnarray*}
and as $x\to 0$ then $w\to 0$.

The problem is now how to choose $\omega$.

From Proposition 2.3 and Proposition 3.1 in Hern\'andez-del-Valle we have 
\begin{eqnarray*}
-\frac{\partial v}{\partial t}+xf''(t)v&=&\frac{1}{2}\frac{\partial^2v}{\partial x^2}+\left(\frac{1}{x}-\frac{x}{s-t}\right)\frac{\partial v}{\partial x}\\
v(s,x)&=&1
\end{eqnarray*}
and
\begin{eqnarray}\label{inequality}
0\leq v(t,x)\leq 1
\end{eqnarray}
admits:
\begin{eqnarray*}
v(t,x):=\tilde{\mathbb{E}}^{t,x}\left[\exp\left\{-\int_t^sf''(u)\tilde{X}_udu\right\}\right].
\end{eqnarray*}
and given that
$$h(t,x)=\frac{x}{\sqrt{2\pi t^3}}\exp\left\{-\frac{x^2}{2t}\right\}$$
then
\begin{eqnarray*}
v(t,x)=\frac{w(t,x)}{h(s-t,x)}
\end{eqnarray*}
where $w$ is an (\ref{pdew1}). From (\ref{inequality})
\begin{eqnarray*}
0\leq w(t,x)\leq h(s-t,x)
\end{eqnarray*}
and hence as $x\to 0$ then $w\to 0$. This will happen for all $t$ if $B_2=0$. Else it holds for $t=0$.

\end{document}